\documentclass{amsart}
\newcommand{\timestamp}{Feburuary 05, 2001}
\usepackage{amsthm}
\usepackage{epic,eepic}
\usepackage[dvips]{graphicx}

\makeatletter

\newcommand{\bmath}[1]{\mbox{\mathversion{bold}$#1$}}

\newcommand{\C}{\bmath{C}}

\newcommand{\Z}{\bmath{Z}}
\newcommand{\R}{\bmath{R}}

\newcommand{\CP}{\bmath{C\!P}}
\newcommand{\SL}{\operatorname{SL}}
\newcommand{\SU}{\operatorname{SU}}
\newcommand{\U}{\operatorname{U}}

\newcommand{\PSU}{\operatorname{PSU}}
\newcommand{\cmcone}{\mbox{\rm CMC}-$1$}
\newcommand{\TA}{\operatorname{TA}}
\newcommand{\trans}{{}^t_{}\!}
\newcommand{\gO}{\mathbf{O}}
\newcommand{\gI}{\mathbf{I}}
\newcommand{\gGamma}{\bmath{\Gamma}}
\newcommand{\ord}{\operatorname{ord}}
\newcommand{\id}{\operatorname{id}}
\newcommand{\Hyp}{\mathcal{H}}
\newcommand{\metone}{\operatorname{Met}_1}
   \newtheorem{theorem}{Theorem}[section]
   \newtheorem*{theorem*}{Theorem}
   \newtheorem{proposition}[theorem]{Proposition}
   \newtheorem{corollary}[theorem]{Corollary}
   \newtheorem{lemma}[theorem]{Lemma}
 \theoremstyle{definition}
   
 \theoremstyle{remark}
   
   \newtheorem*{remark*}{Remark}
\numberwithin{equation}{section}
\makeatother
\title[CMC-1 surfaces of low total curvature II]{
   Mean curvature 1 surfaces in hyperbolic
   3-space with low total curvature II
}
\author{Wayne Rossman}
\author{Masaaki Umehara}
\author{Kotaro Yamada}
\date{\timestamp}
\address[Rossman]{%
   Department of Mathematics, Faculty of Science,
   Kobe University,
   Rokko, Kobe 657-8501, Japan%
}
\email{wayne@math.kobe-u.ac.jp}
\address[Umehara]{%
   Department of Mathematics, Faculty of Science,
   Hiroshima University,
   Higashi-Hiroshima 739-8526, Japan%
}
\email{umehara@math.sci.hiroshima-u.ac.jp}
\address[Yamada]{%
   Faculty of Mathematics,
   Kyushu University 36, 6-10-1
   Hakozaki, Higashi-ku, Fukuoka 812-8185, Japan%
}
\email{kotaro@math.kyushu-u.ac.jp}
\begin{document}
\begin{abstract}
  In this work, complete 
  constant mean curvature $1$ (\cmcone{}) surfaces in hyperbolic
  $3$-space with total absolute curvature at most $4\pi$
  are classified.  This classification suggests that 
  the Cohn-Vossen inequality can be sharpened 
  for surfaces with odd numbers of ends, and a proof of this is given.  
\end{abstract}
\maketitle
\section{Introduction}

This is a continuation (Part~II) of the paper \cite{ruy3} (Part~I) with the 
same title.  
As pointed out in Part~I, complete \cmcone{} (constant mean curvature
$1$) surfaces $f$ in the hyperbolic $3$-space $H^3$ have two important 
invariants.  
One is the {\em total absolute curvature\/} $\TA(f)$, and the 
other is the {\em dual total absolute curvature\/} $\TA(f^\#)$, which is
the total absolute curvature of the dual surface $f^\#$.  
In Part~I, we investigated surfaces with low $\TA(f^\#)$.  
Here we investigate \cmcone{} surfaces with low $\TA(f)$.  

Classifying \cmcone{} surfaces in $H^3$ with low $\TA(f)$ is more
difficult than classifying those with low $\TA(f^\#)$, for the following
reasons: 
$\TA(f)$ equals the area of the spherical image of the (holomorphic)
secondary Gauss map $g$, and $g$ might not be single-valued on the
surface.  
Therefore, $\TA(f)$ is generally not a $4\pi$-multiple of an integer, 
unlike the case of $\TA(f^\#)$.  
Furthermore, the Osserman inequality does not hold for $\TA(f)$, 
also unlike the case of $\TA(f^\#)$.  
The weaker Cohn-Vossen inequality is the best general lower bound for
$\TA(f)$ (with equality never holding \cite{uy1}).  
In Section~\ref{sec:four-pi}, we shall prove the following:  
\begin{theorem}\label{thm:four-pi}
 Let $f:M^2\to H^3$ be a complete \cmcone{} immersion
 of total absolute curvature $\TA(f) \leq 4\pi$.  Then $f$ is either 
 \begin{itemize}
  \item a horosphere,
  \item an Enneper cousin,
  \item an embedded catenoid cousin, 
  \item a finite $\delta$-fold covering of an embedded catenoid cousin
	with $M^2=\C \setminus \{ 0 \}$ and secondary Gauss map
	$g=z^\mu$ for $\mu \leq 1/\delta$, or 
  \item a warped catenoid cousin with injective secondary Gauss map.
 \end{itemize}
\end{theorem}

The horosphere is the only flat (and consequently totally umbilic)
\cmcone{} surface in $H^3$.  
The catenoid cousins are the only \cmcone{} surfaces of 
revolution \cite{Bryant}.  The Enneper cousins are isometric to minimal 
Enneper surfaces \cite{Bryant}.  The warped catenoid cousins \cite{uy1}
are less well known and are described in Section~\ref{sec:prelim}.  

Although this theorem is simply stated, for the reasons stated above the
proof is more delicate than it would be if the condition 
$\TA(f)\le 4\pi$ were replaced with $\TA(f^\#)\le 4\pi$, 
or if minimal surfaces in $\R^3$ with $\TA\le 4\pi$ were considered.  
\cmcone{} surfaces $f$ with $\TA(f^\#)\le 4\pi$ are shown in Part~I to
be only horospheres, Enneper cousin duals, catenoid cousins, and warped
catenoid cousins with embedded ends.  
It is well known that the only complete minimal surfaces in $\R^3$ with
$\TA \leq 4\pi$ are the plane, the Enneper surface, and the catenoid.  

We see from this theorem that any three-ended surface $f$ satisfies 
$\TA(f) > 4\pi$, and so the Cohn-Vossen inequality is not sharp for such
$f$.  
On the other hand, the Cohn-Vossen inequality is sharp for catenoid cousins, 
and a numerical experiment in \cite{ruy4} shows it to be 
sharp for genus $0$ surfaces with $4$ ends.  
This raises the question: 
\begin{quote}
 Which classes of surfaces $f$ have a stronger lower bound for $\TA(f)$ than
 that given by the Cohn-Vossen inequality?
\end{quote}
Pursuing this, in Section~\ref{sec:stronger} we show that stronger 
lower bounds exist for genus zero \cmcone{} surfaces with an 
odd number of ends.

We extend Theorem \ref{thm:four-pi} in a follow-up work \cite{ruy4}, 
to find an inclusive list of 
possibilities for \cmcone{} surfaces with $\TA(f)\le 8\pi$,
and we consider which possibilities we can classify or find examples
for.  
(Minimal surfaces in $\R^3$ with $\TA\le 8\pi$ are classified by Lopez
 \cite{Lopez}.  
 Those with $\TA\le 4\pi$ are listed in Table \ref{tab:minimal}.)
\section{Preliminaries}
\label{sec:prelim}

Let $f\colon{}M\to H^3$ be a conformal \cmcone{} immersion of a Riemann 
surface $M$ into $H^3$.  
Let $ds^2$, $dA$ and $K$ denote the induced metric, 
induced area element and Gaussian curvature, respectively.  
Then $K \leq 0$ and $d\sigma^2:=(-K)\,ds^2$ is a conformal pseudometric
of constant curvature $1$ on $M$.  We call the developing map 
$g\colon{}\widetilde M:=(\mbox{the universal cover of $M$})\to \CP^1$ 
the {\em secondary Gauss map\/} of $f$.  
Namely, $g$ is a conformal map so that its pull-back of the Fubini-Study
metric of $\CP^1$ equals $d\sigma^2$:
\begin{equation}\label{eq:pseudo}
    d\sigma^2 = (-K)\,ds^2 = \frac{4\,dg\,d\bar g}{(1+g\bar g)^2}\;.
\end{equation}
By definition, the secondary Gauss map $g$ of the immersion $f$ is
uniquely determined up to transformations of the form 
\begin{equation}\label{eq:g-change}
   g\mapsto a\star g :=\frac{a_{11}g + a_{12}}{a_{21}g +a_{22}}\qquad
   a=\begin{pmatrix}
        a_{11} & a_{12} \\
        a_{21} & a_{22}
     \end{pmatrix}\in\SU(2)\;.
\end{equation}

In addition to $g$, two other holomorphic invariants $G$ and $Q$ are
closely related to geometric properties of \cmcone{} surfaces.  
The {\em hyperbolic Gauss map\/} $G\colon{}M\to \CP^1$ is holomorphic
and is defined geometrically by identifying the ideal boundary of $H^3$
with $\CP^1$: $G(p)$ is the asymptotic class of the normal geodesic of
$f(M)$ starting at $f(p)$ and oriented in the mean curvature vector's
direction.  
The {\em Hopf differential} $Q$ is the symmetric holomorphic
$2$-differential on $M$ such that $-Q$ is the $(2,0)$-part of the
complexified second fundamental form.
The Gauss equation implies 
\begin{equation}\label{eq:metrics-rel}
  ds^2 \cdot d\sigma^2 = 4 \, Q \cdot \overline Q\;,
\end{equation}
where $\cdot$ means the symmetric product.  
Moreover, these invariants are related by 
\begin{equation}\label{eq:schwarz}
   S(g)-S(G) = 2 Q\;,
\end{equation}
where $S(\cdot)$ denotes the Schwarzian derivative
\[
   S(h):=
    \left[
        \left(\frac{h''}{h'}\right)'-
        \frac{1}{2}\left(\frac{h''}{h'}\right)^2
    \right]\,dz^2\qquad \left ('=\frac{d}{dz}\right)
\]
with respect to a complex coordinate $z$ on $M$.

Since $K\leq 0$, we can define the {\em total absolute curvature} as 
\[
    \TA(f):=\int_M (-K)\,dA\in [0,+\infty] \;\; .
\]
Then $\TA(f)$ is the area of the image in $\CP^1$ of the secondary
Gauss map.  
$\TA(f)$ is generally not an integer multiple of $4\pi$ --- 
for catenoid cousins \cite[Example~2]{Bryant} and their $\delta$-fold 
covers, $\TA(f)$ admits {\em any\/} positive real number.  

For each conformal \cmcone{} immersion $f\colon{}M\to H^3$, there is a
holomorphic null immersion $F\colon{}\widetilde M\to\SL(2,\C)$, 
the {\em lift\/} of $f$, satisfying the differential equation 
\begin{equation}\label{eq:bryant}
   dF = F \begin{pmatrix} 
             g & -g^2 \\ 1 & -g\hphantom{^2}
          \end{pmatrix}\omega \;\; ,\qquad \omega=\frac{Q}{dg}
\end{equation}
such that $f=FF^{*}$, where $F^{*}=\trans\overline{F}$.
Here we consider $H^3=\SL(2,\C)/\SU(2)=\{aa^{*}\,|\,a\in\SL(2,\C)\}$.
If $F=(F_{ij})$, equation \eqref{eq:bryant} implies 
\[ 
   g=-\frac{dF_{12}}{dF_{11}} =-\frac{dF_{22}}{dF_{21}} \; , 
\] 
and it is shown in \cite{Bryant} that 
\[ 
   G = \frac{dF_{11}}{dF_{21}} = \frac{dF_{12}}{dF_{22}} \; . 
\] 

We now assume that the induced metric $ds^2$ on $M$ is complete and 
that $\TA(f)<\infty$, hence there exists a compact Riemann surface 
$\overline M_\gamma$ of genus $\gamma$ and a finite set of points 
$\{p_1,\dots,p_n\}\subset \overline M_\gamma$ ($n\geq 1$) so that 
$M$ is biholomorphic to $\overline M_\gamma\setminus\{p_1,\dots,p_n\}$.
We call the points $p_j$ the {\em ends\/} of $f$.  
Moreover, the pseudometric $d\sigma^2$ as in \eqref{eq:pseudo} is an
element of $\metone(\overline M_{\gamma})$ (\cite[Theorem~4]{Bryant}, 
for a definition of $\metone$ see Appendix~\ref{app:A}).

Unlike the Gauss map for minimal surfaces with $\TA < \infty$ in 
$\R^3$, the hyperbolic Gauss map $G$ of $f$ might not extend to a
meromorphic function on $\overline M_\gamma$ (as the Enneper cousins
show).  
However, the Hopf differential $Q$ does extend to a meromorphic
differential on $\overline M_\gamma$ \cite{Bryant}.  
We say an end $p_j$ $(j=1,\dots,n)$ of a \cmcone{} immersion is 
{\em regular\/} if $G$ is meromorphic at $p_j$.  When $\TA(f) < \infty$,
an end $p_j$ is regular precisely when the order of $Q$ at $p_j$ is at
least $-2$, and otherwise $G$ has an essential singularity at $p_j$
\cite{uy1}. 

Thus the orders of $Q$ at the ends $p_j$ are important for understanding
the geometry of the surface, so we now introduce a notation that
reflects this.  
We say a \cmcone{} surface is of {\em type $\gGamma(d_1,\dots,d_n)$} if
it is given as a conformal immersion 
$f\colon{}\overline M_\gamma\setminus\{p_1,\dots,p_n\}\to H^3$, 
where $\ord_{p_j} Q = d_j$ for $j=1,\dots,n$ 
(for example, if $Q=z^{-2}dz^2$ at $p_1 = 0$, then $d_1=-2$).  
We use $\gGamma$ because it is the capitalized form of $\gamma$, the
genus of $\overline M_\gamma$.  
For instance, $\gI(-4)$ is the class of surfaces of genus $1$ with $1$
end so that $Q$ has an order $4$ pole at the end, 
and $\gO(-2,-3)$ is the class of surfaces of genus $0$ with two ends so
that $Q$ has an order $2$ pole at one end and an order $3$ pole at the
other.
\begin{figure}
\begin{center}
\begin{tabular}{ccc}
       \includegraphics[width=1.3in]{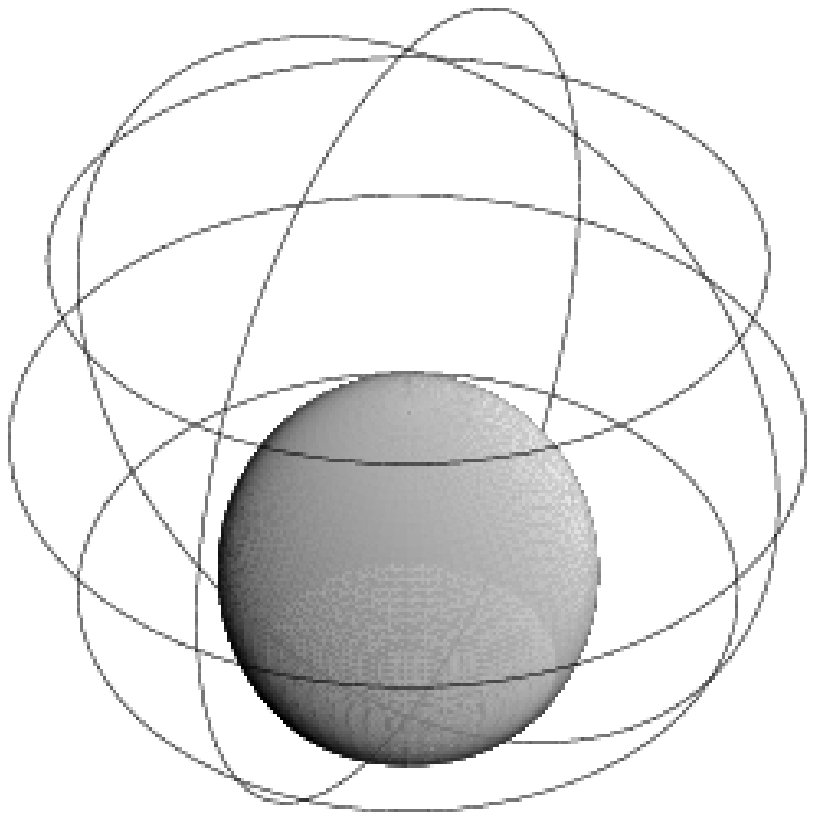}& 
       \includegraphics[width=1.3in]{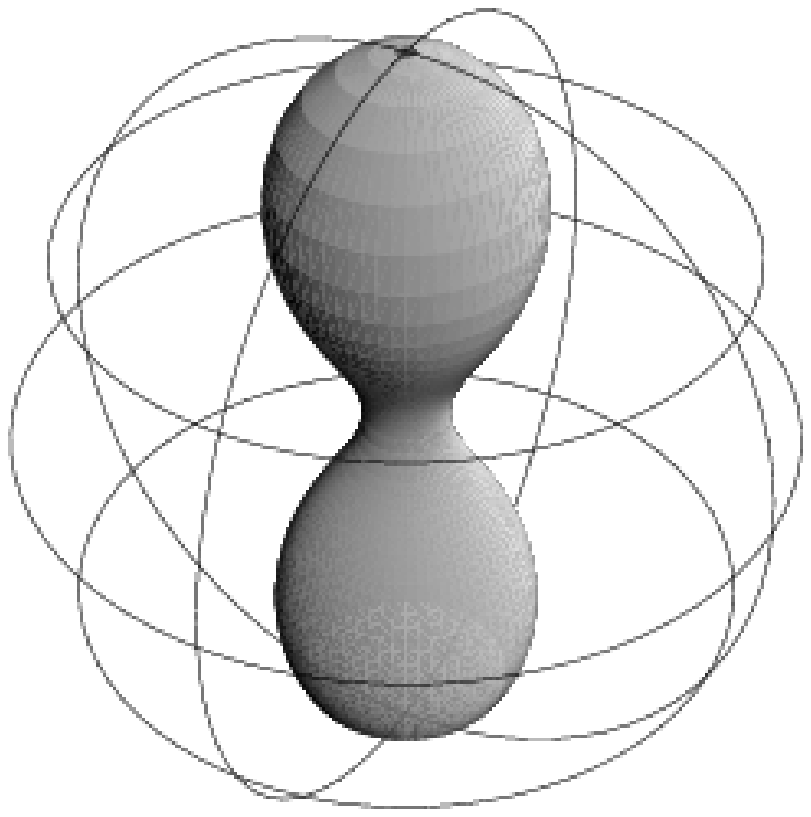} &
       \includegraphics[width=1.3in]{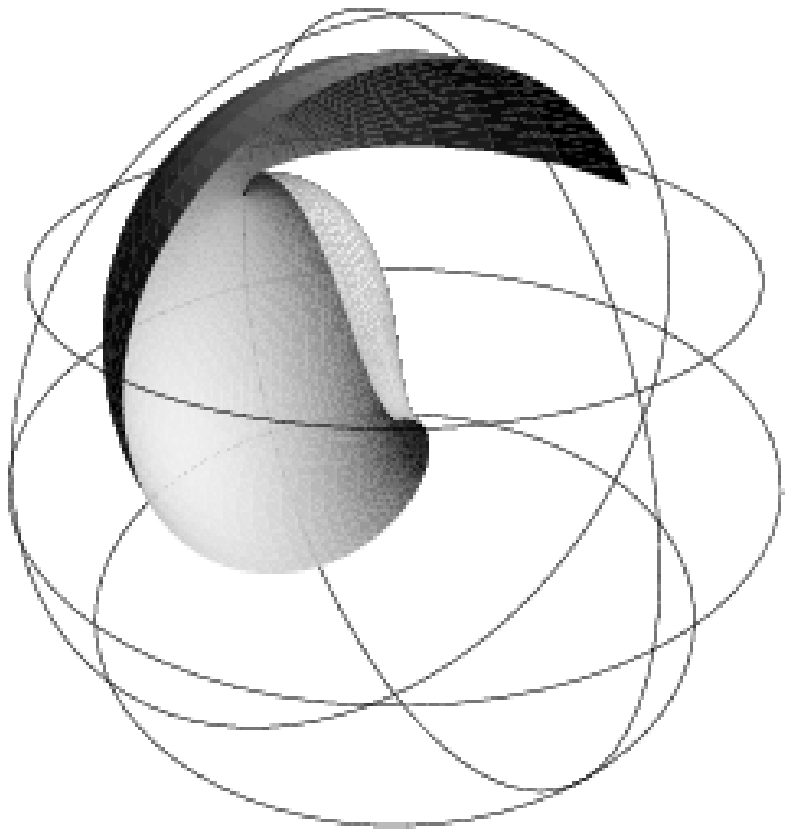}
\end{tabular}
\end{center}
\caption{
   A horosphere, a catenoid cousin with $g=z^\mu$ ($\mu=0.8$), 
   and a fundamental piece 
   (one-fourth of the surface with the end cut away) 
   of an Enneper cousin with $g=z$, $Q=(1/2) dz^2$.
}
\label{fig:enneper}
\end{figure}
\begin{figure}
\begin{center}
\begin{tabular}{cc}
       \includegraphics[width=1.3in]{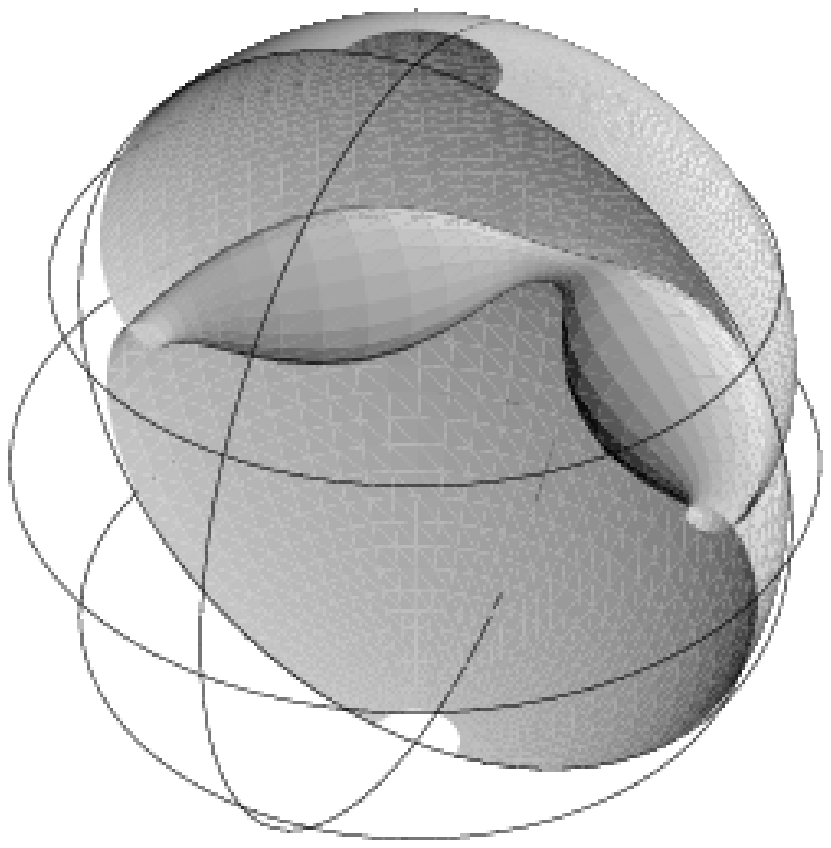} &
       \includegraphics[width=1.3in]{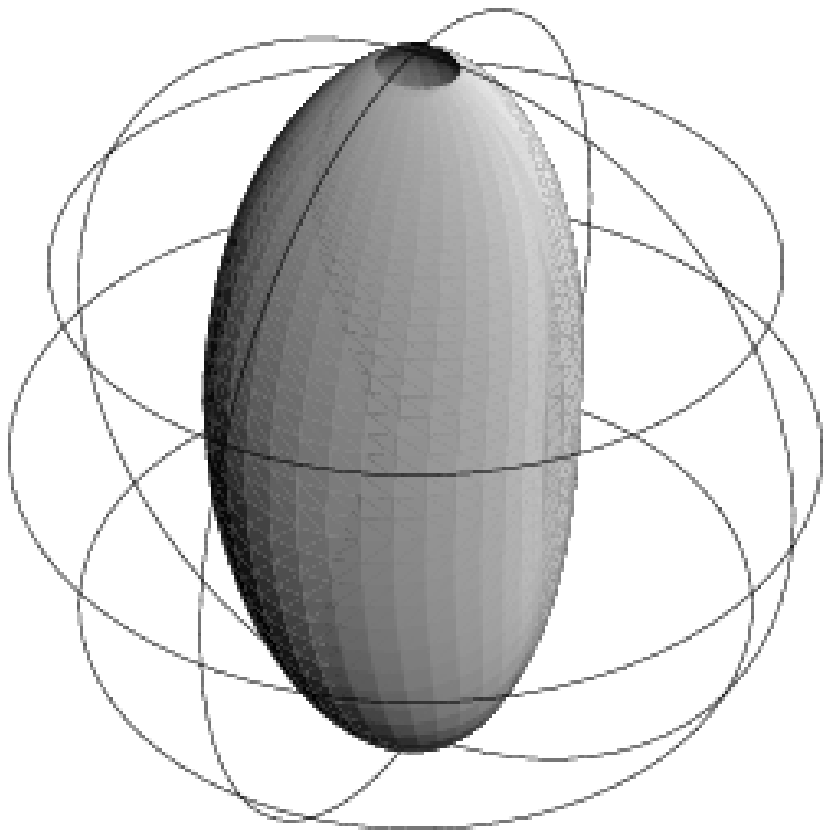}
\end{tabular}
\end{center}
\caption{
 Two warped catenoid cousins, the first with 
 $\delta=1$, $l=4$, $b=1/2$ and the second with $\delta=2$, 
 $l=1$, $b=1/2$.  (Half of the first surface has been cut away.)  
 Only the second of these two surfaces has $\TA(f)=4\pi$ (since $l=1$),
 even though  its ends are not embedded.
}
\label{fig:warped}
\end{figure}
\begin{figure}
  \begin{center}
   \begin{tabular}{ccc}
         \includegraphics[width=1.3in]{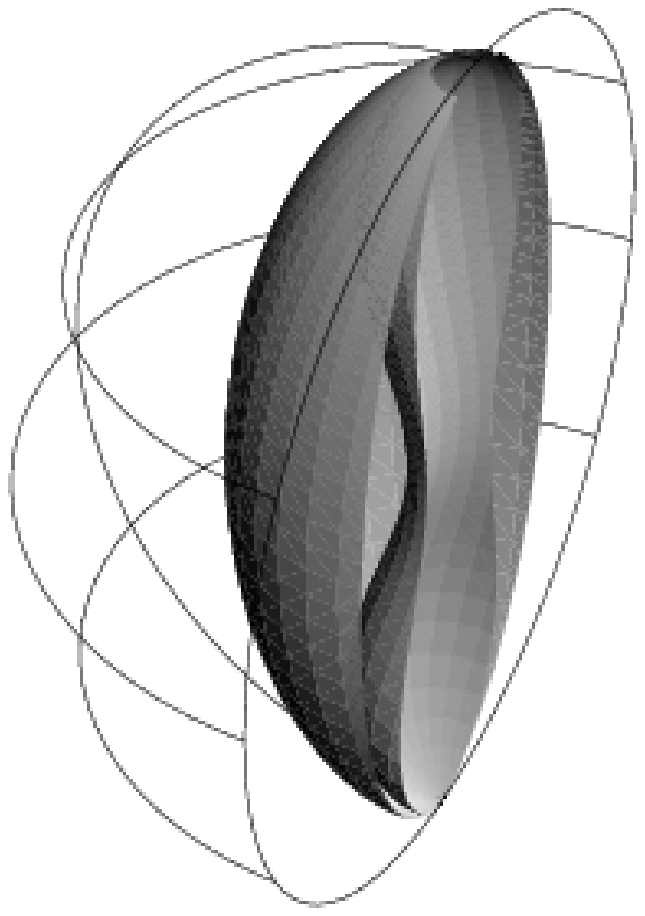} & 
         \includegraphics[width=1.3in]{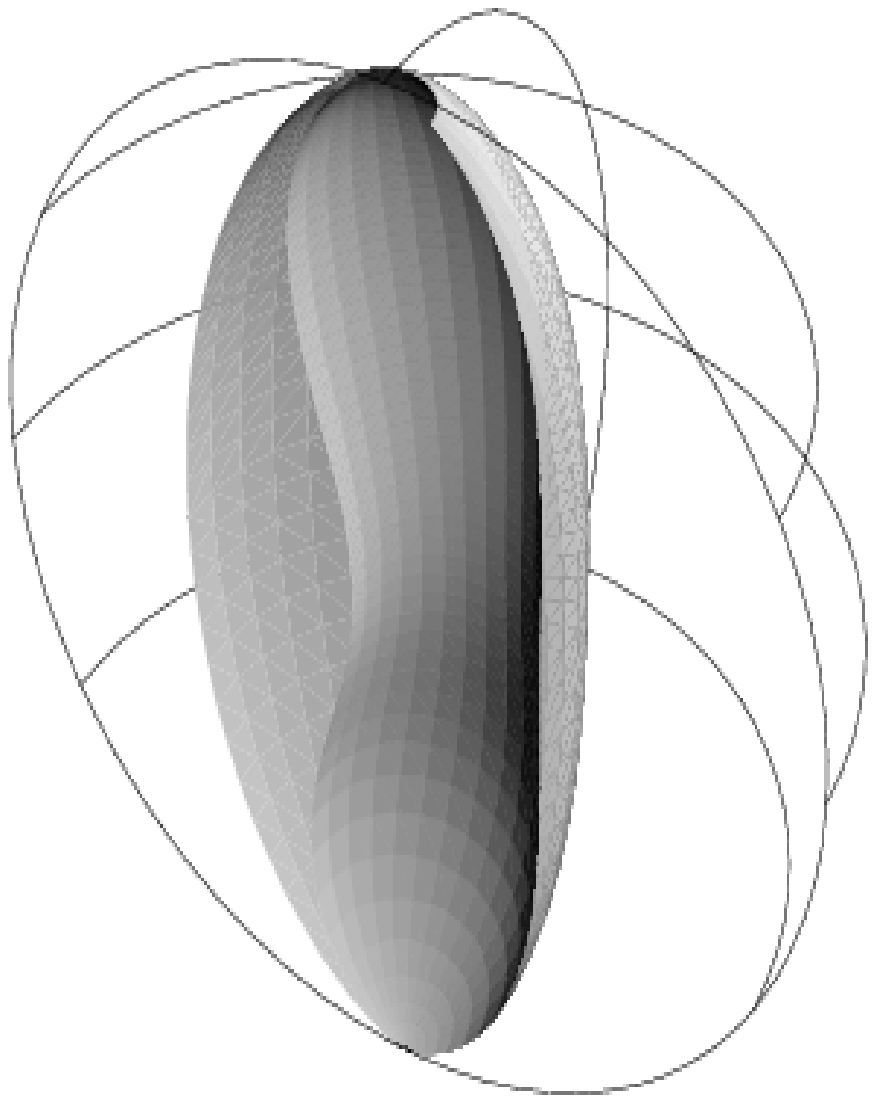} &
         \includegraphics[width=1.3in]{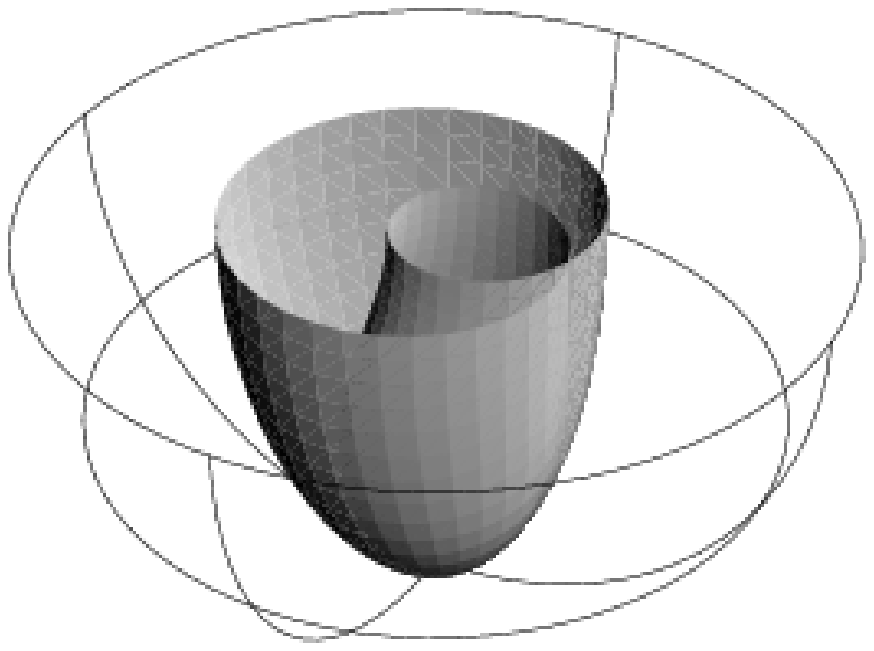}
   \end{tabular}
  \end{center}
\caption{
 Cut-away views of the second warped catenoid cousin in
 Figure~\ref{fig:warped}.
}
\label{fig:warped2}
\end{figure}

We close this section with a description of the warped catenoid cousins.  
Here is a slightly refined version of Theorem~6.2 in \cite{uy1}:  
\begin{theorem}\label{thm:catenoid}
 A complete conformal \cmcone{} immersion $f:M=\C\setminus \{0\} \to H^3$ 
 with two regular ends is a $\delta$-fold cover of a catenoid cousin 
 {\rm(}which is characterized by 
 $g=z^{\mu}$ and $\omega=(1-\mu^2)z^{-\mu-1}dz/(4 \mu)$ for
 $\mu\in\R$\/{\rm)},
 or an immersion {\rm(}or possibly a finite covering of it\/{\rm)}
 where $g$ and $\omega$ can be chosen as 
\[ 
    g = \frac{\delta^2-l^2}{4l}z^l+b \; ,\qquad
    \omega=\frac{Q}{dg} = z^{-l-1}dz \; , 
\] 
 with $l,\delta \in \Z^+$, $l \neq \delta$, and $b \geq 0$.  

 When $b=0$, $f$ is a $\delta$-fold cover of a catenoid cousin with
 $\mu=l$.
 When $b>0$, we call $f$ a {\em warped catenoid cousin}, and its 
 discrete symmetry group is the natural $\Z_2$ extension of the dihedral
 group  $D_l$.  
 Furthermore, the warped catenoid cousins can be written explicitly as 
\[ 
      f= F F^{*},\qquad F = F_0 B \; , 
\]
where 
\[ 
    F_0 = \sqrt{\frac{\delta^2-l^2}{\delta}}
              \begin{pmatrix}
                 \dfrac{1}{l-\delta}    z^{(\delta-l)/2}  &
                 \dfrac{\delta-l}{4l}   z^{(l+\delta)/2} \\
                 \dfrac{1}{l+\delta}    z^{-(l+\delta)/2} &
                 \dfrac{-(l+\delta)}{4l}z^{(l-\delta)/2}
              \end{pmatrix} \quad 
          \text{and} \quad
        B =   \begin{pmatrix}
                    1 & -b \\
                    0 & \hphantom{-}1
              \end{pmatrix}   \,.
\]
\end{theorem}
\begin{proof}
 In \cite{uy1} it is shown that a complete conformal \cmcone{} immersion 
 of $M=\C \setminus \{0\}$ with regular ends is a finite cover of a
 catenoid cousin or an immersion determined by 
\[
 g = az^l+\hat{b} \; , \qquad \omega = c z^{-l-1}dz \; , 
\] 
 where $l$ is a nonzero integer and $a$, $\hat{b}$ and $c$ are complex 
 numbers, which satisfy $l^2+4acl = \delta^2$ for a positive 
 integer $\delta$ and $a,c \neq 0$.  
 (The proof in \cite{uy1} contains typographical errors:
 The exponents $\mu$ and $-\mu$ in equations (6.13) and (6.14) should be
 reversed.  If $\mu \not\in \Z^+$, then the last paragraph of Case 1 is
 correct.  If $\mu \in \Z^+$, then one must consider a possibility that 
 is included in Case 2 in that proof, and the result follows.)  
 Changing $z$ to $1/z$ if necessary, we may assume $l \geq 1$.  

 Choose $\theta$ so that $b:=\hat{b}e^{2i\theta} \geq 0$.  
 Doing the $\SU(2)$ transformation 
\[ 
   g \mapsto 
      \begin{pmatrix}
         e^{i\theta} & 0 \\ 0 & e^{-i\theta} 
      \end{pmatrix} \star g \; , \qquad 
   \omega \mapsto
     e^{-2i\theta} \omega \; , 
\] 
 and replacing $z$ with $e^{-2i\theta/l} c^{1/l} z$ produces the same 
 surface, and one has 
\[
   g = acz^l + b \; , \qquad
   \omega = z^{-l-1}dz \; , \qquad
   ac = \frac{\delta^2 - l^2}{4 l} \; . 
\]  
 Thus $g$ and $\omega$ are as desired.  

 To study the symmetry group of the surface, we consider the 
 transformations 
\[
   \phi_\varrho(z)=e^{2\pi i \varrho/l} \bar{z} \; 
   \quad (\varrho \in \Z) \; ,
   \qquad
   \text{and}
   \qquad 
    \phi(z) = \left( \frac{16l^2(1+b^2)}
     {(\delta^2-l^2)^2}
     \right)^{1/l} \frac{1}{\bar{z}} 
\] 
 of the plane.  
 Then the Hopf differential and secondary Gauss map change as 
\[
       \overline{Q\circ\phi_\varrho}=Q,\qquad
       \overline{g\circ\phi_\varrho}= g,\qquad
       \overline{Q\circ\phi}=Q,\qquad
       \overline{g\circ\phi}= \frac{bg+1}{g-b}=A\star g\;,
\] 
 where 
\[
       A = \frac{i}{\sqrt{1+b^2}}
           \begin{pmatrix}
                b & \hphantom{-}1 \\
                1 & -b
            \end{pmatrix} \in \SU(2)\;.
\]
 Hence $\phi_{\varrho}$ and $\phi$ represent isometries of the surface.  
 One can then check that there are no other isometries of the surface, 
 i.e., that there are no other anti-conformal bijections $\hat{\phi}$ of
 $M$ so that  $\overline{Q \circ \hat{\phi}} = Q$ and 
 $\overline{g \circ \hat{\phi}} =  A\star g$ for some $A \in SU(2)$.  
 Thus the symmetry group is $D_l \times  \Z_2$.  

 To see that the warped catenoid cousins have the explicit
 representation  described in the theorem, one needs only to verify that
 $F=F_0 B$ satisfies \eqref{eq:bryant}.  
\end{proof}
\begin{table}\rm
  \begin{tabular}{|l||c|l|}
    \hline
      Type & $\TA$ & The surface \\
    \hline\hline
      $\gO(0)$  & 0 
                & Plane \\
    \hline
      $\gO(-4)$ & $4\pi$ 
                & Enneper surface \\
    \hline
      $\gO(-2,-2)$ 
                & $4\pi$
                & Catenoid \\
    \hline
  \end{tabular}
\caption{Classification of minimal surfaces in $\R^3$ with 
$\TA\leq 4\pi$.}
\label{tab:minimal}
\end{table}

\section{Complete CMC-$1$ surfaces with $\TA(f) \leq 4\pi$}
\label{sec:four-pi}
\begin{table}
\begin{center}
  \begin{tabular}{|l||c|l|}
    \hline
     Type & $\TA(f)$ & The surface \\
    \hline\hline
      $\gO(0)$  & 0 
                & Horosphere \\
    \hline
      $\gO(-4)$ & $4\pi$ 
                & Enneper cousins   \\
    \hline
      $\gO(-2,-2)$ 
                & $(0,4\pi]$
                & Catenoid cousins and     \\
                & 
                & \hspace{2em}their $\delta$-fold covers   \\
    \hline
      $\gO(-2,-2)$ 
                & $4\pi$
                & Warped catenoid  \\
                & 
                & \hspace{2em}cousins with $l=1$   \\
    \hline
  \end{tabular}
\end{center}
\caption{Classification of CMC-$1$ surfaces in $H^3$ with 
$\TA(f) \leq 4\pi$.}
\label{tab:ta-4pi}
\end{table}
In this section we will prove Theorem~\ref{thm:four-pi}.  
First we prepare some notations and basic facts.  
For a complete conformal \cmcone{} immersion 
$f\colon{}M=\overline M_{\gamma}\setminus\{p_1,\dots,p_n\}\to H^3$, we
define $\mu_j$ and $\mu^{\#}_j$ to be the branching orders of the Gauss
maps $g$ and $G$, respectively, at each end $p_j$.  
At an irregular end $p_j$, we have $\mu_j^{\#}=\infty$.  
At a regular end $p_j$, relation \eqref{eq:schwarz} implies that
the Hopf differential $Q$ expands as 
\begin{equation}\label{eq:q-first}
    Q = \left(\frac{1}{2}\frac{c_j}{(z-p_j)^2} + \dots \right)\,dz^2 \; , 
    \qquad
    c_j = -\frac{1}{2}\mu_j(\mu_j+2) +
           \frac{1}{2}\mu_j^{\#}(\mu_j^\#+2) \; , 
\end{equation}
where $z$ is a local complex coordinate around $p_j$.

Let $\{q_1,\dots,q_m\} \subset M$ be the $m$ umbilic points of the surface, 
and let $\xi_k=\ord_{q_k}Q $.  Then, as in (2.5) of Part~I, 
\begin{equation}\label{eq:fact-a}
   \sum_{j=1}^n d_j+\sum_{k=1}^m 
             \xi_k=4\gamma-4, 
   \qquad
   \text{in particular,}
   \quad  \sum_{j=1}^n d_j\le 4\gamma-4\;.
\end{equation}
By \eqref{eq:metrics-rel} and \eqref{eq:schwarz}, it holds that
\begin{equation}\label{eq:fact-h}
    \xi_k = \text{[branch order of $G$ at $q_k$]}
           = \text{[branch order of $g$ at $q_k$]}
          = \ord_{q_k} d\sigma^2 \, . 
\end{equation}
As in (2.4) of Part~I, the Gauss-Bonnet theorem implies 
\[
   \frac{\TA(f)}{2\pi}=\chi(\overline M_{\gamma})+
                    \sum_{j=1}^{n}\mu_j + \sum_{k=1}^m \xi_k \; ,
\]
where $\chi$ denotes the Euler characteristic.  
Combining this with \eqref{eq:fact-a}, we have 
\begin{equation}\label{eq:fact-b}
  \frac{\TA(f)}{2\pi}=2\gamma-2+\sum_{j=1}^n(\mu_j-d_j)\;.
\end{equation}
Proposition~4.1 in \cite{uy1} implies that
\begin{equation}\label{eq:fact-c}
   \mu_j-d_j>1, \qquad
   \text{in particular,} \qquad
   \mu_j-d_j\geq 2\quad 
   \text{if $\mu_j\in\Z$}\;.
\end{equation}
An end $p_j$ is regular if and only if $d_j\geq -2$, and then $G$ is 
meromorphic at $p_j$. Thus 
\begin{equation}\label{eq:fact-f}
  \text{$\mu_j^{\#}$ is a non-negative integer if $d_j\geq -2$}.
\end{equation}
By Proposition 4 of \cite{Bryant}, 
\begin{equation}\label{eq:fact-g}
   \mu_j > -1 \; , 
\end{equation}
hence equation \eqref{eq:q-first} implies 
\begin{equation}\label{eq:fact-d}
   \mu_j=\mu_j^\#\in\Z\qquad 
   \text{if $d_j\geq -1$}\;.
\end{equation}
Finally, we note that
\begin{equation}\label{eq:fact-e}
  \parbox[t]{10cm}{
    any meromorphic function on a Riemann surface 
    $\overline M_{\gamma}$ of genus $\gamma\geq 1$
    has at least three distinct branch points.
  }
\end{equation}
To prove this, let $\varphi$ be a meromorphic function on 
$\overline M_\gamma$ with $N$ branch points $\{q_1, \dots, q_N\}$ of
branching order $\psi_k$ at $q_k$.  
Then the Riemann-Hurwicz relation implies 
\[
   2\deg \varphi = 2 - 2\gamma + \sum_{k=1}^{N}{\psi_k}\;.
\]
On the other hand, since the multiplicity of $\varphi$ at $q_k$ 
is $\psi_k+1$, $\deg\varphi\geq \psi_k+1$ ($k=1,\dots,N$).  Thus 
\[
   (N-2)\deg\varphi \geq 2(\gamma-1)+N\;.
\]
If $\gamma\geq 1$, then $\deg \varphi \geq 2$, and so $N \geq 3$.  
\begin{remark*}
 Facts \eqref{eq:fact-b} and \eqref{eq:fact-c} imply that, 
 for \cmcone{} surfaces, equality never holds in the Cohn-Vossen 
 inequality \cite{uy1}:
 \begin{equation}\label{eq:oss}
  \frac{\TA(f)}{2\pi} > -\chi(M)=n-2+2\gamma \; . 
 \end{equation}
\end{remark*}

\begin{proof}[Proof of Theorem~\rm{\ref{thm:four-pi}}]
By \eqref{eq:fact-b},
\begin{equation}\label{eq:4}
  2\ge \frac{\TA(f)}{2\pi}=2\gamma-2+\sum_{j=1}^{n}(\mu_j-d_j)\;.
\end{equation}
Since $\mu_j-d_j>1$ by \eqref{eq:fact-c}, we have
\[
   4> 2\gamma + n\;.
\]
Thus the only possibilities are 
\[
  (\gamma,n)=(0,1),\quad (0,2),\quad (0,3),\quad (1,1)\;.
\]

\subsection*{The case \boldmath$(\gamma,n)=(1,1)$}
By \eqref{eq:4} and \eqref{eq:fact-g}, we have $d_1\ge \mu_1-2>-3$.
Thus the end $p_1$ is regular, and $G$ is meromorphic on $\overline
M_1$.  By \eqref{eq:fact-a}, $d_1 \le 0$. 
If $d_1=-2$, the end has non-vanishing flux, 
and the surface does not exist, by Corollary~3 of \cite{ruy2}. 
If $d_1=0$ or $-1$, by \eqref{eq:fact-a} there is at most one umbilic
point.  
Since any branch point of $G$ is at an end or an umbilic point, 
\eqref{eq:fact-e} is contradicted.  Hence a surface of this type does
not exist.  

\subsection*{The case \boldmath$(\gamma,n)=(0,1)$}
Here the surface is simply connected, so it has a canonical
isometrically corresponding minimal surface in $\R^3$ with the same
total absolute curvature.  
We conclude the surface is a horosphere or an Enneper cousin.

\subsection*{The case \boldmath$(\gamma,n)=(0,2)$}
Here, by \eqref{eq:fact-a}, we have $d_1+d_2\le -4$.  
On the other hand, by \eqref{eq:4} and \eqref{eq:fact-g}, we have
$d_1+d_2\ge -4+(\mu_1+\mu_2)>-6$.  
Thus $d_1+d_2$ is either $-4$ or $-5$, and we now consider these two
cases separately:  
\subsubsection*{The case \boldmath$d_1+d_2=-4$}
If $d_1+d_2=-4$, then there are no umbilic points, by \eqref{eq:fact-a}.
If $d_1,d_2\ge -2$, then the ends are regular, and Theorem~\ref{thm:catenoid} 
implies the surface is a $\delta$-fold cover of an embedded catenoid cousin 
with $\delta \leq 1/\mu$, or a warped catenoid cousin with $l=1$.  

Now assume that
\[
   d_1\ge -1,\qquad d_2\le -3.
\]
Then we have $\mu_1\in \Z$ by \eqref{eq:fact-d}.  
By Proposition~\ref{prop:a-1} in Appendix~\ref{app:A},
we cannot have just one $\mu_j \not\in \Z$, so also $\mu_2\in \Z$. 
Then $g$ is single-valued on $M$.  Since $g$ and $G$ are both
single-valued on $M$, the lift $F$ is also 
(see equations (1.6) and (1.7) in \cite{uy3}), 
and so the dual immersion $f^{\#}$ is also single-valued on $M$.  
Since $(f^{\#})^\#=f$, $f^{\#}$ is a \cmcone{} immersion with dual total
absolute curvature $4\pi$ and of type $\gO(-1,-3)$.  Such an $f^\#$
cannot exist by Theorem~3.1 of Part~I, so such an $f$ does not exist.  

\subsubsection*{The case \boldmath$d_1+d_2=-5$}
If $d_1+d_2=-5$, then the surface has only one umbilic point $q_1$ with 
$\xi_1 = 1$, by \eqref{eq:fact-a}, and we can set $\overline M_0=
\C\cup\{\infty\}$, $p_1=0$, $p_2=\infty$, and $q_1=1$.  

By \eqref{eq:4}, $\mu_1+\mu_2\leq -1$. 
Then by \eqref{eq:fact-g} at least one of $\mu_1$ and $\mu_2$ is not an 
integer.  
Hence both are not integers, by Proposition~\ref{prop:a-1} in
Appendix~\ref{app:A}.  
Then \eqref{eq:fact-d} implies we may assume $d_1=-2$ and $d_2=-3$.  
by Proposition~\ref{prop:a-2} in Appendix~\ref{app:A}, 
the metric $d\sigma^2$ is the pull-back of the Fubini-Study metric on
$\CP^1$ by the map 
\[
  g=c z^\mu\left(z-\frac{\mu+1}{\mu}\right) \qquad 
  \left(
    c\in \C\setminus\{0\},\, \mu\in\R\setminus\{0,\pm 1\}
  \right)\;.
\]
On the other hand, the Hopf differential is of the form
\begin{equation}\label{eq:o23q}
  Q(z)=q(z)\,dz^2=
  \theta\frac{z-1}{z^2}\,dz^2 \qquad (\theta\in \C\setminus\{0\}).
\end{equation}
Thus $\omega=Q/dg$ can be written in the form
\begin{equation}\label{eq:o23omega}
   \omega=w(z)\,dz=
       \frac{\theta}{c}\frac{1}{\mu+1}\frac1{z^{\mu+1}}\,dz\;.
\end{equation}
Consider the equation (which is introduced in \cite{uy1} as (E.1))
\begin{equation}\label{eq:e-1}
  X''+a(z)X'+b(z)X=0\;,\qquad
   \left(
    a(z):=-\frac{w'(z)}{w(z)},\quad
    b(z):=-q(z)
   \right)\;.
\end{equation}
We expand the coefficients $a$ and $b$ as
\[
  a(z)=\frac{1}{z}\sum_{j=0}^\infty a_j z^j,\qquad
        b(z)=\frac{1}{z^2}\sum_{j=0}^\infty b_j z^j\;,
\]
Then the origin $z=0$ is a regular singularity of equation \eqref{eq:e-1}.
Let $\lambda$ and $\lambda+m$ be the solutions of the corresponding
indicial equation $t(t-1)+a_0 t + b_0=0$ with $m \geq 0$.  
If the surface exists, then Theorem~2.4 of \cite{uy1} implies 
that $m$ must be a positive integer and the log-term coefficient 
of the solutions of \eqref{eq:e-1} must vanish.  
When $m \in \Z^+$, the log-term coefficient vanishes if and only if 
\[
        \sum_{k=0}^{m-1}
                 \left\{(\lambda+k)a_{m-k}+b_{m-k}\right\}
                 \eta_k
                 (\lambda) = 0\;,
\]
where $\eta_0 = 1$ and $\eta_1,\dots,\eta_{m-1}$ are given recursively by 
\[
       \eta_j =\frac{1}{j(m-j)}\displaystyle\sum_{k=0}^{j-1}
           \left\{(\lambda+k)a_{j-k}+b_{j-k}\right\}
           \eta_k
\]
as in Proposition~A.3 in Appendix~A of Part~I.  Here we have
\[
   0=a_1=a_2=\ldots \; , \qquad 0=b_2=b_3=\ldots \; , 
\]
and so the log-term coefficient 
never vanishes at the end $p_1$, because $b_1=-\theta 
\neq 0$.  Thus this type of surface does not exist.  

\subsection*{The case \boldmath$(\gamma,n)=(0,3)$}
This is the only remaining case.  But this type of surface 
does not exist, by the following Theorem \ref{thm:threenoid}. 
\end{proof}
\begin{theorem}\label{thm:threenoid}
 Let $f\colon{}M\to H^3$ be a complete \cmcone{} immersion of genus zero
 with three ends. Then $\TA(f) > 4\pi$.  
\end{theorem}

\begin{remark*}
   The second and third authors proved that $\TA(f)\geq 4\pi$
   holds for \cmcone{} surfaces of genus $0$ with three ends
   \cite[Proposition~2.7]{uy6}.
   Then the essential part of Theorem~\ref{thm:threenoid} is that
   $\TA(f)=4\pi$ is impossible.
\end{remark*}

\begin{proof}[Proof of Theorem~\ref{thm:threenoid}]
We suppose $\TA(f) = 4\pi$, and will arrive at a contradiction. 
Without loss of generality, we may set $\overline M_0=\C\cup\{\infty\}$
and $p_1=0$, $p_2=1$, and $p_3=\infty$.  
\subsection*{Step 1}
Since $\gamma=0$ and $\TA(f) \leq 4\pi$, \eqref{eq:fact-b} implies 
\begin{equation}\label{3}
   4\ge \sum_{j=1}^3(\mu_j-d_j) \; . 
\end{equation}
Since $\mu_j-d_j >1$ for all $j$, \eqref{3} implies $\mu_j-d_j <2$ for all 
$j$.  
Hence $\mu_1,\mu_2,\mu_3\not\in \Z$ by \eqref{eq:fact-c}.  
Then \eqref{eq:fact-d} implies $d_j \le -2$ for all $j$, and as 
equations \eqref{3} and \eqref{eq:fact-g} imply 
$d_1+d_2+d_3 \ge -4+\mu_1+\mu_2+\mu_3 > -7$, 
we have 
\begin{equation}\label{6}
   d_1=d_2=d_3=-2 \; , 
\end{equation}
and so the ends are regular.  

On the other hand, since $\TA(f)=4\pi$, \eqref{eq:fact-b} and \eqref{6}
imply 
\begin{equation}\label{7}
    \mu_1+\mu_2+\mu_3=-2\;.
\end{equation}
Then by \eqref{eq:fact-g}, we have
\begin{equation}\label{7.5}
   -1<\mu_j<0\qquad (j=1,2,3) \; , 
\end{equation}
and furthermore at least two of the $\mu_j$ are less than $-1/2$.  
We may arrange the ends so that 
\begin{equation}\label{eq:murange}
  -1<\mu_1,\mu_2 < -\frac{1}{2} \quad
   \text{and}
  \quad
  -1<\mu_3 < 0 \; . 
\end{equation}
Moreover, by Appendix~A of \cite{uy6} (note that the $C_j$ there equal 
$\pi (\mu_j+1)$), the metric $d\sigma^2$ is reducible (as defined in 
Appendix~\ref{app:B} of the present paper).  
Then, by Proposition~\ref{prop:B1} and the relation \eqref{eq:order-diff} 
in the appendices here, the secondary Gauss map 
$g$ can be expressed in the form 
\begin{equation}\label{g}
   g=z^{-(\mu_1+1)} (z-1)^{\beta+1} \frac{a(z)}{b(z)} \; , 
\end{equation}
where $a(z)$, $b(z)$ are relatively prime polynomials without 
zeros at $p_1$ and $p_2$, and 
\begin{equation}\label{beta}
  \beta=\mu_2 \qquad \text{or} \qquad \beta=-2-\mu_2 \; .
\end{equation}
Note that the order of $g$ at $p_3=\infty$ is $\pm (\mu_3+1)$ and is
also $\mu_1-\beta-\deg a + \deg b$.
If $\beta=\mu_2$, 
\[ 
  2\mu_1 = \deg a -\deg b - 1\qquad \text{or}\qquad
  2\mu_2 = \deg b -\det a - 1
\]
holds.
Thus either $2\mu_1$ or $2\mu_2$ is an integer,
but this contradicts \eqref{eq:murange}, so $\beta = -\mu_2 - 2$:
\begin{equation}\label{eq:first-g}
    g = z^{-\mu_1-1}(z-1)^{-\mu_2-1}\frac{a(z)}{b(z)}\;.
\end{equation}
Thus, by \eqref{7}, we have
\[ 
   -\mu_3-\deg a+\deg b= \pm (\mu_3+1) \; . 
\]  
Hence either 
\begin{gather}
 \deg a-\deg b=1 \qquad 
  \text{and the order of $g$ at $\infty$ is $-\mu_3-1$, or}
\label{eq:first-case}
\\
 \mu_3=-1/2, \deg a=\deg b
     \qquad
   \text{and the order of $g$ at $\infty$ is $\mu_3+1$}
\label{eq:second-case}
\end{gather}
holds because of \eqref{eq:murange}.
To get more specific information about $a(z)$ and $b(z)$, we now consider 
$dg$: 

\subsection*{Step 2}
Since $Q$ is holomorphic on $\C\setminus\{0,1\}$ with 
two zeroes (by \eqref{eq:fact-a}), \eqref{eq:q-first} implies 
\begin{equation}\label{eq:Q}
  Q=\frac{1}{2}\left(
               \frac{c_3 z^2+(c_2-c_1-c_3)z+c_1}{z^2(z-1)^2}
            \right)\,dz^2\;,
\end{equation}
with the $c_j$ as in \eqref{eq:q-first}, as pointed out in 
\cite[page~84]{uy6}.
Note that 
\begin{equation}\label{eq:positive-cj}
  c_j > 0 \qquad (j=1,2,3) \; , 
\end{equation}
because $\mu_j^\#\geq 0$ and $-1<\mu_j<0$.  
Let $q_1$ and $q_2$ be the two roots 
of 
\begin{equation}\label{quadraticequation}
   c_3 z^2+(c_2-c_1-c_3)z+c_1=0\;.
\end{equation} 
In the case of a double root, we write $q:=q_1=q_2$.

Using \eqref{eq:fact-h} and Proposition~\ref{prop:B1} in
Appendix~\ref{app:B}, $dg$ has only the four following 
possibilities:  
\begin{align}
  dg&= C\frac{z^{-\mu_1-2} (z-1)^{-\mu_2-2} (z-q_1)(z-q_2) }{
          \prod_{k=1}^r(z-a_k)^2}\,dz\;, 
\label{eq:dg-possible-1}\\
  dg&= C \frac{z^{-\mu_1-2}(z-1)^{-\mu_2-2} (z-q_1)}{
              (z-q_2)^3\prod_{k=1}^r(z-a_k)^2}\,dz
               \qquad &(q_1\neq q_2)\;,
\label{eq:dg-possible-2}\\
  dg&=C \frac{z^{-\mu_1-2}(z-1)^{-\mu_2-2}}{
             (z-q_1)^3 (z-q_2)^3\prod_{k=1}^r(z-a_k)^2}\,dz
      \qquad &(q_1\neq q_2)\;,
\label{eq:dg-possible-3}
\intertext{or}
  dg &= C \frac{z^{-\mu_1-2} (z-1)^{-\mu_2-2}}{
           (z-q)^4\prod_{k=1}^r(z-a_k)^2}\,dz
        \qquad &(q=q_1=q_2)\;,
\label{eq:dg-possible-4}
\end{align}
where $r$ is a non-negative integer and the points 
$a_k \in \C \setminus \{0,1,q_1,q_2\}$ are mutually distinct.  
In the first case \eqref{eq:dg-possible-1}, 
the order of $dg$ at infinity ($z=p_3=\infty$) is given by
\[
   \mu_1+\mu_2+2r=2r-2-\mu_3=\mu_3\text{ or }-\mu_3-2 \; . 
\]
So $2r-2=2 \mu_3 \in(-2,0) $ or $2r-2-\mu_3=-\mu_3-2$.
Hence $r=0$ and the order of $dg$ at $\infty$ is
$-\mu_3-2$ in the first case.

In the other three cases \eqref{eq:dg-possible-2},
\eqref{eq:dg-possible-3}, \eqref{eq:dg-possible-4}, the orders of $dg$
at infinity are 
\[
   \mu_1+\mu_2+(2\text{ or }6\text{ or }4)+2r+2 \geq 2-\mu_3+2r > 2 \; , 
\]
respectively.  
These orders must equal either $\mu_3<0$ or $-\mu_3-2<0$, so none of these 
three cases can occur.  We conclude that $dg$ is of the form
\begin{equation}\label{dg}
   dg= C z^{-\mu_1-2} (z-1)^{-\mu_2-2} (z-q_1)(z-q_2)\, dz
     \qquad (C\in\C\setminus\{0\})\;.
\end{equation}
Since the order of $dg$ at $\infty$ is $\mu_1+\mu_2=-\mu_3-2<0$,
\eqref{eq:first-case} holds.

\subsection*{Step 3}
Now we determine the polynomials $a(z)$, $b(z)$ in expression
\eqref{eq:first-g}. 
Differentiating \eqref{eq:first-g}, we have
\begin{equation}\label{ex:dg}
   dg=\frac{z^{-\mu_1-2}(z-1)^{-\mu_2-2}}{b^2(z)} f(z) dz \;,
\end{equation}
where 
\begin{equation}\label{f}
  f(z)=-(1+\mu_1)(z-1)ab-(1+\mu_2)z ab+z(z-1)(a'b-ab')\;.
\end{equation}
Since $a(z)$ and $b(z)$ are relatively prime, $b(z)$ does not divide 
$f(z)$ when $\deg b \ge 1$. But \eqref{dg} and \eqref{ex:dg} imply that 
$b^2(z)$ divides $f(z)$, so $b(z)$ is constant, and we may assume $b=1$.  
Here, as seen in the previous step, \eqref{eq:first-case} holds, and then,
$\deg a = 1$.
Thus we have
\begin{equation}\label{eq:ab}
   a(z)= a_1 z + a_0 \qquad \text{and}\qquad b=1 \qquad (a_1\neq 0)\;.
\end{equation}
\subsection*{Step 4:}
By \eqref{dg}, \eqref{ex:dg}, \eqref{f} and \eqref{eq:ab} we have 
\begin{multline}\label{identity}
-a_1(\mu_1+\mu_2+1)z^2+\left\{\mu_1a_1-(\mu_1+\mu_2+2)a_0\right\}z
      +(1+\mu_1)a_0\\
  = C (z-q_1)(z-q_2) \; . 
\end{multline}
Equation \eqref{quadraticequation} also has roots $q_1$ and $q_2$, so 
\begin{equation}\label{eq:sol-coef}
   q_1q_2=\frac{a_0}{a_1}\frac{1+\mu_1}{1+\mu_3}=
         \frac{c_1}{c_3},\quad
   q_1+q_2=
        -\frac{\mu_3a_0+\mu_1a_1}{a_1(1+\mu_3)}
           =\frac{c_1}{c_3}+1-\frac{c_2}{c_3}\;.
\end{equation}
By \eqref{eq:fact-g}, \eqref{eq:positive-cj} and the first equation of 
\eqref{eq:sol-coef}, we have $a_0/a_1>0$.  
Substituting the first equation of \eqref{eq:sol-coef} into the second, 
we have
\[
   \frac{c_2}{c_3}=-\frac{1+\mu_2}{1+\mu_3}\left(
              \frac{a_0}{a_1}+1\right).
\]
Since $a_0/a_1>0$, \eqref{eq:fact-g} implies $c_2/c_3<0$, contradicting 
\eqref{eq:positive-cj} and completing the proof.  
\end{proof}
\section{Improvement of the Cohn-Vossen Inequality}
\label{sec:stronger}

For a complete \cmcone{} immersion $f$ into $H^3$,
equality in the Cohn-Vossen inequality never holds
(\cite[Theorem~4.3]{uy1}).  In particular, when $f$ is of 
genus $0$ with $n$ ends, 
\begin{equation}\label{eq:cohn-vossen}
    \TA(f)>2\pi (n-2)\;.
\end{equation}
For $n=2$, the catenoid cousins show 
that \eqref{eq:cohn-vossen} is sharp.  But Theorem \ref{thm:threenoid} 
implies 
\[
   \TA(f)>4\pi \qquad\text{for}\qquad n=3 \; , 
\]
which is stronger than the Cohn-Vossen inequality \eqref{eq:cohn-vossen}.
The following theorem gives a sharper inequality than that of 
Cohn-Vossen when $n$ is any odd integer:  
\begin{theorem}\label{thm:odd}
  Let $f\colon{}\C\cup\{\infty\}\setminus\{p_1,\dots,p_{2l+1}\}\to H^3$
  be a complete conformal genus $0$ \cmcone{} immersion with $2l+1$ ends, 
  $l \in \Z$.  Then 
\[
  \TA(f)\ge 4\pi l\;.
\]
\end{theorem}

To prove this, we first prepare two lemmas and a proposition. 
\begin{lemma}\label{lem:trig}
  Let $\theta_1,\theta_2,\theta_3\in [0,\pi]$ 
  be three real numbers such that 
\begin{equation}\label{th1}
  \cos^2\theta_1+\cos^2\theta_2+\cos^2\theta_3+
     2\cos\theta_1\cos\theta_2\cos\theta_3\le 1 \; . 
\end{equation}
 Then the following inequalities hold\/{\rm : }
\begin{align}
  \theta_1+\theta_2+\theta_3&\ge \pi \; , \label{eq:t1} \\
  \theta_2-\theta_1 &\le \pi-\theta_3 \; . \label{eq:t2} 
\end{align}
\end{lemma}
\begin{remark*}
  It is well known that the inequality
\[
   \cos^2\theta_1+\cos^2\theta_2+\cos^2\theta_3+
       2\cos\theta_1\cos\theta_2\cos\theta_3< 1
\]
  is a necessary and sufficient condition for the 
  existence  of a spherical triangle $\mathcal T$ with 
  angles $\theta_1$, $\theta_2$ and $\theta_3$.
  Then \eqref{eq:t1} follows directly 
  from the Gauss-Bonnet formula, and 
  \eqref{eq:t2} is the triangle inequality
  for the polar triangle of $\mathcal T$, and the lemma follows.  
  ($\mathcal T$'s polar triangle is the one whose vertices are the centers of 
  the great circles containing the edges of $\mathcal T$.)  
  However, we give an alternative proof: 
\end{remark*}
\begin{proof}[Proof of Lemma~{\rm\ref{lem:trig}}]
 We set
\[
  E:=\cos^2\theta_1+\cos^2\theta_2+\cos^2\theta_3+
       2\cos\theta_1\cos\theta_2\cos\theta_3- 1 \leq 0 \; . 
\]
 Then 
\begin{multline*}
    E =   4 \cos\left(\frac{\theta_1+\theta_2+\theta_3}2\right)
           \cos\left(\frac{-\theta_1+\theta_2+\theta_3}2\right) \times \\
       \cos\left(\frac{\theta_1-\theta_2+\theta_3}2\right)
       \cos\left(\frac{\theta_1+\theta_2-\theta_3}2\right) \; . 
\end{multline*}
 If $\theta_1+\theta_2+\theta_3<\pi$, we have
$|\pm\theta_1\pm\theta_2\pm\theta_3|<\pi$, and so 
\[
   \cos\left(\frac{\pm\theta_1\pm\theta_2\pm\theta_3}2\right)> 0 \; , 
\]
 implying $E>0$, a contradiction.  This proves \eqref{eq:t1}.  Now, since 
\[
   E=\cos^2\theta_1+\cos^2(\pi-\theta_2)+\cos^2(\pi-\theta_3)
    +2\cos\theta_1\cos(\pi-\theta_2)\cos(\pi-\theta_3)- 1
\]
 and $E\le 0$ and $\theta_1,\pi-\theta_2,\pi-\theta_3\in [0,\pi]$, 
  \eqref{eq:t1} implies 
\[
   \theta_1+(\pi-\theta_2)+(\pi-\theta_3)\ge \pi \; , 
\]
 that is, \eqref{eq:t2} holds.  
\end{proof}

For a matrix $a\in\SU(2)$, there is a unique 
$C\in[0,\pi]$ so that $a$ has eigenvalues $\{-e^{\pm iC}\}$.
We define the {\em rotation angle\/} of $a$ as 
\[
    \theta(a):=2C\;.
\]
In fact, if one considers the matrix acting on the unit sphere as an 
isometry (M\"obius action on $\CP^1$ with the Fubini-Study metric), 
$\theta(a)$ is exactly the angle of rotation.
\begin{lemma}\label{lem:key}
 Let $a_0$, $a_1$, $a_2$, $a_3$ be four matrices in $\SU(2)$ satisfying
 $a_1a_2a_3=a_0$. Then it holds that 
\[
  \theta(a_1)+\theta(a_2)+\theta(a_3)\ge \theta(a_0) \; . 
\]
\end{lemma}
\begin{proof}
 Setting $b:=a_3(a_0)^{-1}=(a_1a_2)^{-1}$, we have $a_1a_2b=\id$.  
 Then Appendix~A of \cite{uy6} implies 
\[
   \cos^2 \frac{\theta(a_1)}2+
     \cos^2 \frac{\theta(a_2)}2+\cos^2 \frac{\theta(b)}2
            +2\cos \frac{\theta(a_1)}2
      \cos\frac{\theta(a_2)}2\cos\frac{\theta(b)}2\le 1 \; . 
\]
 So by Lemma \ref{lem:trig} we have
\begin{equation}\label{eq:angle1}
   \frac{\theta(a_1)}2+
   \frac{\theta(a_2)}2+
   \frac{\theta(b)}2\ge \pi \; . 
\end{equation}
 On the other hand, we have $a_3^{-1}ba_0=\id$.
 Again Appendix A of \cite{uy6} implies 
\[
   \cos^2 \frac{\theta(a_0)}2+
      \cos^2 \frac{\theta(a_3)}2+
      \cos^2 \frac{\theta(b)}2+
      2\cos \frac{\theta(a_0)}2
        \cos \frac{\theta(a_3)}2\cos \frac{\theta(b)}2\le 1 \; , 
\]
 since $\theta(a_3^{-1})=\theta(a_3)$.  
 By \eqref{eq:t2} of Lemma \ref{lem:trig}, we have
\begin{equation}\label{eq:angle2}
  \frac{\theta(a_0)}{2}-\frac{\theta(a_3)}{2}
       \leq \pi-\frac{\theta(b)}{2} \; . 
\end{equation}
 By \eqref{eq:angle1} and \eqref{eq:angle2}, we
 get the assertion.
\end{proof}

\begin{proposition}\label{prop:key}
 Let $a_1,\dots,a_{2m+1}$ be matrices in $\SU(2)$ satisfying 
\[
      a_1a_2\cdots a_{2m+1}=\id\;. 
\]
  Then it holds that 
\[
   \sum_{j=1}^{2m+1} \theta(a_j)\ge 2\pi \; . 
\]
\end{proposition}

\begin{remark*}  This result does not hold for an even number of matrices: 
 Suppose $a_1$, \dots, $a_{2m}$ $\in \SU(2)$ satisfy 
 $a_1a_2\cdots a_{2m}=\id$.  Then the inequality 
 $\sum_{j=1}^{2m} \theta(a_j)\ge 0$ is sharp. In fact, equality will hold if 
 all $a_j=-\id$.
\end{remark*}
\begin{proof}[Proof of Proposition~{\rm\ref{prop:key}}]
 We argue by induction.  
 If $m=1$, the result follows from Lemma~\ref{lem:key} with $a_0=\id$.  
 Now suppose the result always holds for $m-1(\ge 1)$.  Set 
\[
   b:=a_1a_2a_3 \; . 
\]
 Then, by Lemma \ref{lem:key}, 
\begin{equation}\label{theta1}
 \theta(a_1)+\theta(a_2)+\theta(a_3)\ge \theta(b) \; .
\end{equation}
 On the other hand, we have $b a_4\cdots a_{2m+1}=\id$, so 
by the inductive assumption, 
\begin{equation}\label{theta2}
 \theta(b)+\sum_{j=4}^{2m+1} \theta(a_j)\ge 2\pi \; . 
\end{equation}
By \eqref{theta1} and \eqref{theta2}, we get the assertion.  
\end{proof}
We now apply Proposition~\ref{prop:key} to the monodromy representation 
of pseudometrics in $\metone(\C\cup\{\infty\})$ (see Appendices~\ref{app:A} 
and \ref{app:B}): 
\begin{corollary}\label{cor:odd}
 Let $d\sigma^2 \in \metone(\C\cup\{\infty\})$ with divisor 
\[
    D= \sum_{j=1}^{s} \beta_j p_j 
           + \sum_{k=1}^{n} \xi_k q_k \; , \;\;\;\;\; \beta_j>-1 \; , \;\;\; 
              \xi_k \in \Z^+ \; , 
\]
 where the $p_1,\dots,p_{s},q_1,\dots q_n$ are mutually distinct
 points in $\C\cup\{\infty\}$.

 If $s+\xi_1 + \dots + \xi_n$ is an odd integer, then 
$\beta_1 + \dots + \beta_{s} \geq 1-s$.  
\end{corollary}
\begin{proof}
 Let $g$ be a developing map of $d\sigma^2$ with  the monodromy representation 
 $\rho_g\colon{}\pi_1(M)\to\PSU(2)=\SU(2)/\{\pm\id\}$ on 
 $M=\C\cup\{\infty\}\setminus\{p_1,\dots,p_s,q_1,\dots,q_n\}$.
 
 $\rho_g$ can be lifted to an $\SU(2)$ representation 
 $\tilde\rho_g\colon{}\pi_1(M)\to\SU(2)$ so that the two following properties 
 hold: 
 \begin{enumerate}
 \item \label{property1} Let $T_j$ ($j=1,\dots,s$) and $S_k$ 
 ($k=1,\dots,n$) be deck transformations on $\widetilde M$ corresponding 
 to loops about $p_j$ and $q_k$, respectively.  Then it holds that 
\[
    \tilde{\rho}_{g}(T_1)\dots\tilde{\rho}_{g}(T_s)
    \tilde{\rho}_{g}(S_1)\dots\tilde{\rho}_{g}(S_n)=\id \; . 
\]
 \item The eigenvalues of the matrix $\tilde \rho_g(T_j)$ 
 (resp. $\tilde \rho_g(S_k)$)
 are $\{-e^{\pm i \pi(\beta_j+1)}\}$ (resp. $\{-e^{\pm i\pi(\xi_k+1)}\}$).
 \end{enumerate}
 
 This is proven in \cite[Lemma~2.2]{uy6} for $s=3,n=0$, and the same 
 argument will work for general $s$ and $n$.  We include an outline of the 
 argument here:  one chooses a solution $\tilde F$ to equation 
 (2.12) in \cite{uy6} (with $G=z$ and $Q=S(g)/2$), then $\tilde F$ 
 has a monodromy representation $\rho_{\tilde F}\colon{}\pi_1(M)\to\SU(2)$, 
 where $\tilde F \to \tilde F \cdot \rho_{\tilde F}(\gamma)$ about loops 
 $\gamma \in \pi_1(M)$.  Then $\rho_g= \pm \rho_{\tilde F}$, and we simply 
 choose the lift $\tilde \rho_g$ so that $\tilde \rho_g= + \rho_{\tilde F}$.  
 The first property is then clear.  
 
 To show the second property, we note that when $\beta_j$ and $\xi_k$ are all 
 given the value $0$, then $Q$ is identically $0$ and so $\tilde F$ is 
 constant and all $\rho_{\tilde F}= + \id$.  Hence the eigenvalues 
 $\{\pm e^{\pm i \pi(\beta_j+1)}\}$ (resp.~$\{\pm e^{\pm i\pi(\xi_k+1)}\}$) 
 of $\tilde \rho_g(T_j)$ (resp.~$\tilde \rho_g(S_k)$) 
 are $\{-e^{\pm i \pi(\beta_j+1)}\}$ (resp.~$\{-e^{\pm i\pi(\xi_k+1)}\}$) 
 in this case.  Then, as $\beta_j$ and $\xi_k$ are deformed 
 back to their original values, the matrices $\tilde \rho_g(T_j)$ 
 (resp.~$\tilde \rho_g(S_k)$) change analytically and so the sign of the 
 eigenvalues cannot change, showing the second property.  
 
 We have
\[
    \theta(\tilde \rho_{g}(T_j)) \leq 2\pi(\beta_j+1)\;,
\]
 and since $\xi_k$ is an integer, we have 
\begin{equation}\label{eq:monodromy-int}
    \tilde{\rho}_{g}(S_k) = (-1)^{\xi_k}\id \; . 
\end{equation}

 Assume $s=2m+1$ is an odd number.
 Then, by the assumption, $\xi_1+\dots+\xi_n$ is an even integer, and 
 by \eqref{eq:monodromy-int} above we have 
 $\tilde{\rho}_{g}(T_1)\dots\tilde{\rho}_{g}(T_{2m+1})=\id$, so by
 Proposition~\ref{prop:key}, 
\[
   2\pi\sum_{j=1}^{2m+1}(\beta_j+1) \geq 
   \sum_{j=1}^{2m+1}\theta(\tilde{\rho}_{g}(T_j)) \geq 2\pi \; , 
\]
 implying the corollary when $s$ is odd.  

 Now suppose that $s=2m$ is even.  We have 
 $\tilde{\rho}_g(S_1)\dots\tilde{\rho}_g(S_n)=-\id$, 
 because $\xi_1+\dots+\xi_n$ is odd. 
 Hence 
 $\tilde{\rho}_g(T_1)\dots\tilde{\rho}_g(T_{2m})
 (-\id) = \id$, and since $\theta(-\id)=0$, Proposition~\ref{prop:key} 
 implies 
\[   2 \pi \sum_{j=1}^{2m} (\beta_j+1) \geq 
    \sum_{j=1}^{2m} \theta(\tilde \rho_g(T_j)) + \theta(-\id) \geq 2\pi \; , 
\]
 implying the corollary when $s$ is even.  
\end{proof}
\begin{proof}[Proof of Theorem~{\rm\ref{thm:odd}}]
 Suppose that $\mu_1 \in \Z$.  Then by \eqref{eq:fact-b} and 
 \eqref{eq:fact-c}, 
\begin{equation}\label{extralabelinproof}
  \frac{\TA(f)}{2\pi}\ge
        -2+ (\mu_1-d_1)+
           \sum_{j=2}^{2m+1}({\mu_j-d_j}) >
            -2 +2+2m=2 m \; , 
\end{equation}
 proving the theorem when $\mu_1 \in \Z$.  
 
 Next, suppose that $d_1\leq -3$.
 In this case, $\mu_1-d_1>-1+3=2$.
 Hence again by \eqref{eq:fact-b} and \eqref{eq:fact-c}, we have
\eqref{extralabelinproof}, and the theorem follows.  

 Thus we may assume $\mu_j \not\in \Z$  and $d_j\geq -2$ at all ends.  
 Then by \eqref{eq:fact-d} we have all $d_j= -2$.
 So, by \eqref{eq:fact-a} and \eqref{eq:fact-h}, the
 corresponding pseudometric $d\sigma^2$ has divisor
\[
    \sum_{j=1}^{2m+1}\mu_jp_j + \sum_{k=1}^{l} \xi_k q_k,\qquad
    \sum_{k=1}^l \xi_k =4m-2\in 2\Z\;,
\]
 where $\xi_k=\ord_{q_k} Q$ at each umbilic point $q_k$
 ($k=1,\dots,l$).  Then by Corollary~\ref{cor:odd},
\[
    \mu_1 + \mu_2 +\dots +\mu_{2m+1}\geq -2m \; , 
\]
and so \eqref{eq:fact-b} implies the theorem.  
\end{proof}
\begin{remark*}
 When $m=1$, we know the lower bound $4\pi m$ in Theorem~\ref{thm:odd}
 is sharp.  However, we do not know if it is sharp for general  $m$.  
 For genus $0$ \cmcone{} surfaces with an even number $n \geq 4$ of ends, 
 we do not know if there exists any stronger lower bound than that of the 
 Cohn-Vossen inequality.  

 In \cite{ruy4}, it is shown numerically that there exist \cmcone{} 
 surfaces of genus $0$ with four ends whose total absolute curvature 
 gets arbitrarily close to $4\pi$.
\end{remark*}
\appendix
\begingroup
\renewcommand{\thesection}{\Alph{section}}
\section{}
\label{app:A}
\endgroup

For a compact Riemann surface $\overline M$ and points $p_1,\dots,p_n \in 
\overline M$, a conformal metric $d\sigma^2$ of constant curvature $1$ 
on $M:=\overline M\setminus\{p_1,\dots,p_n\}$ is an element of 
$\metone(\overline M)$ if there exist real numbers $\beta_1,\dots,\beta_n 
>-1$ so that each $p_j$ is a conical singularity of order $\beta_j$, that is, 
if $d\sigma^2$ is asymptotic to $c_j|z-p_j|^{2\beta_j}\,dz\cdot d\bar z$ 
at $p_j$, for $c_j\neq 0$ and $z$ a local complex coordinate around $p_j$.  
We call the formal sum 
\begin{equation}\label{eq:divisor}
   D:=\sum_{j=1}^n \beta_j p_j 
\end{equation}
the {\em divisor\/} corresponding to $d\sigma^2$.
For a pseudometric $d\sigma^2\in\metone(\overline M)$ with divisor 
$D$, there is a holomorphic map 
$g\colon{}\widetilde M\to\CP^1$ such that $d\sigma^2$ is 
the pull-back of the Fubini-Study metric of 
$\CP^1$.  This map, called the {\em developing map\/} of $d\sigma^2$, 
is uniquely determined up to M\"obius transformations $g \mapsto a\star g$ for 
$a \in \SU(2)$.  

For a conical singularity $p_j$ of $d\sigma^2$, there exists 
a developing map $g$ and a local coordinate $z$ of $\overline M$
around $p_j$ such that 
\[
   g(z) = (z-p_j)^{\tau_j} \hat g(z) \qquad (\tau_j\in\R\setminus\{0\}) \; , 
\]
where $\hat g(z)$ is holomorphic in a neighborhood 
of $p_j$ and $\hat g(p_j)\neq 0$.
Here, the order $\beta_j$ of $d\sigma^2$ at $p_j$ is
\begin{equation}\label{eq:order}
      \beta_j = \begin{cases}
                \hphantom{-}\tau_j-1\qquad &   \text{if $\tau_j>0$}\;, \\
                 -\tau_j-1 \qquad & \text{if $\tau_j<0$}\;. 
                \end{cases}
\end{equation}
In other words, if $dg = (z-p_j)^{\beta} \hat h(z)\,dz$, 
where $\hat h(z)$ is holomorphic near $p_j$ and $\hat h(p_j)\neq 0$, 
then the order $\beta_j$ is expressed as 
\begin{equation}\label{eq:order-diff}
      \beta_j = \begin{cases}
                   \hphantom{-}\beta &   \text{if $\beta>-1$}\;, \\
                   -\beta-2 \qquad &     \text{if $\beta<-1$}\;. 
                \end{cases}
\end{equation}

The following proposition gives an obstruction to the existence of certain 
pseudometrics in $\metone(\C\cup\{\infty\})$.  
\begin{proposition}\label{prop:a-1}
 For any non-integer $\beta>-1$, there is no pseudometric 
 $d\sigma^2$ in $\metone(\C\cup \{\infty\})$ with the divisor 
\[
    \beta p_1+\sum_{j=2}^n m_jp_j\qquad
       (m_2,\dots,m_n\in\Z) \;,
\]
 where $p_1,\ldots,p_n$ are mutually distinct points in 
 $\C\cup\{\infty\}$.
\end{proposition}
When $n=1$ (i.e. when $\sum_{j=2}^n m_jp_j$ is removed), this nonexistence 
 of a ``tear-drop" has been pointed out in \cite{Troyanov1} and \cite{CL}.  
\begin{proof}
  We may set $p_1=\infty$.  Since the $m_j \in \Z$, 
  the developing map $g$ of $d\sigma^2$ is well-defined on 
  $\C$, and so $g$ is meromorphic on $\C$.  
  As $d\sigma^2$ has finite total curvature, 
  $g$ extends to $z=\infty$ as a holomorphic mapping.  
  In particular, $\beta \in \Z$.  
\end{proof}
\begin{remark*}
  When a Riemann surface $\overline{M}_{\gamma}$ has genus $\gamma>0$, 
  there is a pseudometric in $\metone(\overline{M}_\gamma)$ with only 
  one singularity that has order less than $0$, by \cite{Troyanov2}.  
\end{remark*}
\begin{proposition}\label{prop:a-2}
 Suppose a pseudometric $d\sigma^2$ in 
 $\metone(\C\cup \{\infty\})$ has divisor 
 \[
    \beta_1 p_1 + \beta_2 p_2 + p_3 \qquad 
   (\text{$\beta_1,\beta_2>-1$ and $\beta_1,\beta_2\not\in \Z$})
   \;, 
 \]
 where $p_1:=0$, $p_2:=\infty$, and $p_3:=1$. 
 Then $d\sigma^2$ has a developing map $g$ of the form 
\begin{equation}\label{eq:norm-g}
    g=cz^\mu\left(z-\frac{\mu+1}{\mu}\right) \qquad 
       (c\in \C,\, \mu\in \R) \; , 
\end{equation}
 where $\beta_1=|\mu|-1$ and $\beta_2=|\mu+1|-1$.
\end{proposition}
\begin{proof}
  Since $d\sigma^2$ has only two non-integral conical singularities, it is 
  reducible, and 
  Proposition~\ref{prop:B1} in Appendix~\ref{app:B} shows that the map
  $g$ is written in the form
  \[
      g=z^\mu \frac{a(z)}{b(z)} \qquad (\mu\not\in \Z)\;,
  \]
  where $a(z)$ and $b(z)$ are relatively prime polynomials with 
  $a(0)\ne 0$ and $b(0)\ne 0$.  
  Note that $b(z)$ can have a multiple root only at a conical singularity of 
  $d\sigma^2$, hence only at $z=1$.  
  Thus $b'(z_0)\neq 0$ for all roots $z_0\in\C\setminus\{0,1\}$ of $b$.
  
  Since the change $g\mapsto 1/g$ preserves $d\sigma^2$,
  we may assume that $\deg a \geq \deg b$.
  By a direct calculation, we have
  \[
    dg(z)=\frac{z^{\mu-1}}{b(z)^2}h(z) dz \; ,
  \text{   with } h(z)=\mu a(z)b(z)+za'(z)b(z)-za(z)b'(z) \; .
  \]
  Note that $h(0)=\mu a(0) b(0)\neq 0$.
  
  Let $z_0\in\C\setminus\{0,1\}$.
  If $b(z_0) \neq 0$,  then
  $g(z_0) \neq \infty$, and since $z_0$ is not a singularity 
  of $d\sigma^2$, 
  we have $dg(z_0) \neq 0$, and hence $h(z_0) \neq 0$.  
  If $b(z_0) = 0$, then $a(z_0) \neq 0$ and 
  $b'(z_0) \neq 0$, so $h(z_0) \neq 0$.  
  Hence the only root of the polynomial $h(z)$ is $1$:
  \[
     h(z)=k\, (z-1)^m \; , \;\;\; m \in \Z^+ \; , \;\; k \in \C\setminus\{0\} \; . 
  \]
  
  We claim that $m=1$.  
  If $b(1) \neq 0$, then $g$ (or $d\sigma^2$) having 
  order $1$ at $p_3=1$ means that $m=1$, 
  by \eqref{eq:order-diff} and the above form of $dg(z)$.  
  Suppose $b(1)=0$.  Then we have 
  $b(z)=(z-1)^l\hat b(z)$, where $\hat b(z)$ is a polynomial 
  in $z$ with $\hat b(1)\neq 0$ and $l \in \Z^+$.  
  Furthermore, $h(z)=(z-1)^{l-1}\hat h(z)$, where $\hat h(z)$ is
  a polynomial with $\hat h(1)\neq 0$,
  since $a(1)\neq 0$.  So $m=l-1$. 
  Then by \eqref{eq:order-diff}, we have $m=1$.
  
  Suppose that $\deg b\ge 1$.  
  Since $\deg a\ge \deg b$, the top term of $h(z)$ must vanish. 
  Thus we have $\mu=\deg b - \deg a \in \Z$, 
  contradicting that $\beta_1, \beta_2 \not\in \Z$.
  So $b(z)$ is constant.  
  Similarly, if $\deg a \geq 2$, then $\mu = - \deg a \in \Z$.  
  Hence $\deg a = 1$,   and $g$ is as in \eqref{eq:norm-g}.  
  $\beta_1=|\mu|-1$ and $\beta_2=|\mu+1|-1$ follow from 
  \eqref{eq:order-diff}.  
\end{proof}
\begingroup
\renewcommand{\thesection}{\Alph{section}}
\section{}
\label{app:B}
\endgroup
Consider $d\sigma^2\in\metone(\overline M)$ with divisor $D$ as in 
\eqref{eq:divisor} in Appendix~\ref{app:A} and developing map $g$.  
Since the Fubini-Study metric of $\CP^1$ is invariant 
under the deck transformation group $\pi_1(M)$ of 
$M:=\overline M\setminus\{p_1,\dots,p_n\}$, there is 
a representation 
\[
    \rho_g\colon{}\pi_1(M)\longrightarrow 
    \SU(2) 
\]
such that 
\[ 
    g\circ T^{-1}  = \rho_g(T)\star g \qquad (T\in\pi_1(M))\;.
\]
The metric $d\sigma^2$ is called {\em reducible\/} if the image 
of $\rho_g$ is a commutative subgroup in $\SU(2)$, and is called 
{\em irreducible\/} otherwise.  
Since the maximal abelian subgroup of $\SU(2)$ is $\U(1)$, 
the image of $\rho_g$ for a reducible $d\sigma^2$ lies in a subgroup 
conjugate to $\U(1)$, and this image might be simply the identity.  
We call a reducible metric $d\sigma^2$ {\em $\Hyp^3$-reducible\/}
if the image of $\rho_g$ is the identity, and {\em $\Hyp^1$-reducible\/}
otherwise (for more on this, see \cite[Section~3]{ruy1}).

Let $p_1,\dots,p_{n-1}$ be distinct points in $\C$ and $p_n=\infty$.  
We set 
\[
   M_{p_1,\dots,p_n}:=\C\cup\{\infty\} 
     \setminus \{p_1,p_2,\dots,p_n\} \qquad (p_n=\infty)
\]
and $\widetilde{M}_{p_1,\dots,p_n}$ its universal cover.  

The following assertion was needed in the proof of Theorem~\ref{thm:four-pi}. 
\begin{proposition}\label{prop:B1}
  Let $p_1,\dots,p_{n-1}$ be mutually distinct points of $\C$, 
  and let $d\sigma^2$ be a metric of constant curvature $1$ 
  defined on $M_{p_1,\dots,p_n}$ 
  $(p_n=\infty)$ which has a conical singularity at each $p_j$. 
  Suppose $d\sigma^2$ is reducible and 
  $\beta_j:=\ord_{p_j} d\sigma^2$ satisfy 
\[
     \beta_1,\dots,\beta_m\not\in\Z \, ,\qquad
     \beta_{m+1},\dots,\beta_{n-1}\in\Z \, ,\qquad
          \beta_n\not\in\Z \, , 
\]
  for some $m\leq n-1$.  Then the metric $d\sigma^2$ has a developing map 
  $g\colon{}\widetilde M_{p_1,\dots,p_n}\to \C\cup \{\infty\}$ 
  given by 
\[
     g=(z-p_1)^{\tau_1}\cdots (z-p_m)^{\tau_m} r(z)
       \qquad (\tau_1,\dots,\tau_m\in \R\setminus\Z) \; , 
\]
  where $r(z)$ is a rational function on $\C\cup\{\infty\}$ and
\[
    (z-p_1)^{\tau_1}\cdots (z-p_m)^{\tau_m}:=\exp\left(
         \sum_{j=1}^m \tau_j \int_{z_0}^z \frac{dz}{z-p_j}
                   \right) \qquad (z\in M_{p_1,\dots,p_n}),
\]
  for some base point $z_0\in M_{p_1,....,p_n}$.
\end{proposition}
\begin{proof}
  $d\sigma^2$ is reducible only if the image of 
  the representation $\rho_g$ is simultaneously diagonalizable, 
  so we may choose a developing map 
  $g\colon{}\widetilde{M}_{p_1,\dots,p_n} \to \CP^1$ such that
\begin{equation}\label{eq:red-repr}
    \rho_g(T) = 
    \begin{pmatrix}
        e^{i\theta_T} & 0 \\ 
        0 & e^{-i\theta_T} 
    \end{pmatrix}\;.
\end{equation}  
  Thus we have 
\[
    \log(g \circ T^{-1})=\log(g) + 2i\theta_T \; .  
\]
  Differentiating this gives 
\[
   d\log(g \circ T^{-1})=d\log(g) \; , 
\]
  which implies that $d\log(g)$ is single-valued 
  on $M_{p_1,\dots,p_n}$.  

  On the other hand, by Proposition~4 in \cite{Bryant}, there is a
  complex coordinate $w$ around each end $p_j$ such that
\begin{equation}\label{eq:g-norm}
  a_j\star g=(w-p_j)^{\tau_j} \qquad (\tau_j\in \R\setminus\{0, \pm 1 \})\; ,
\end{equation} 
  for some $a_j \in \SU(2)$ ($j=1,\dots,n$).  
  Let $T_j$ be the deck transformation of 
  $\widetilde{M}_{p_1,\dots,p_n}$ corresponding to a loop 
  surrounding $p_j$.  Then 
\[
   \rho_g(T_j)\neq \pm\id\quad \text{ for }j=1,\dots,m\text{ and }j=n \; . 
\]  
  Hence $\tau_j \not\in \Z$ when $j \leq m$ and $j=n$.  
  By \eqref{eq:red-repr}, $a_j$ in \eqref{eq:g-norm} is diagonal, so 
\[
    g(p_j) = 0 \quad \text{or}\quad \infty\qquad
        (j=1,\dots,m,n) \; . 
\]
  Hence $d\log(g)$ has poles of order $1$ at
  $p_1,\dots,p_m$, and thus 
\[
  d\log(g) = 
    \frac{dg}{g} 
     = \frac{\tau_1\, dz}{z-p_1}+\dots+\frac{\tau_m\, dz}{z-p_m}
      +u(z)\, dz\;,
\]
  where $u(z)$ is meromorphic.  Integrating this gives the assertion.  
\end{proof}


\begin{thebibliography}{RUY}
\bibitem{BC}
    J.~L.~M.~Barbosa and A.~G.~Colares,
    {\sc Minimal Surfaces in $\R^3$},
    Lecture Notes in Math. vol.~1195, 
    Springer-Verlag,  (1986).
\bibitem{Bie}
    L.~Bieberbach, 
    {\sc Theorie der gew\"ohnlichen Differential\-gleich\-ungen},
    Zweite Auflage,
    Springer-Verlag, (1965).
\bibitem{Bryant}
  R.~Bryant,
  {\itshape Surfaces of mean curvature one in hyperbolic space},
  Ast\'erisque {\bfseries 154--155} (1987), 321--347.
\bibitem{CL}
  W. Chen and C. Li,
  {\it What kinds of singular surfaces can admit constant curvature?},
  Duke Math.~J., {\bf 78} (1995) 437--451.
\bibitem{cg}
  C.~C.~Chen, F.~Gackstatter,
  {\itshape Elliptische und hyperelliptische Funktionen und
    vollst\"andige Minimalfl\"achen vom Enneperschen Typ},
    Math.~Ann. {\bfseries 259} (1982), 359--369.
\bibitem{chr}
  P.~Collin, L.~Hauswirth and H.~Rosenberg,
  {\itshape The geometry of finite topology Bryant surfaces},
  to appear in Ann. of Math.
\bibitem{ET1} R. Sa Earp, E. Toubiana, {\itshape On the geometry of 
   constant mean curvature one surfaces in hyperbolic space}, to appear in 
   Illinois J. Math.  
\bibitem{ET2} R. Sa Earp, E. Toubiana, {\itshape Meromorphic data for 
   mean curvature one surfaces in hyperbolic space}, preprint.  
\bibitem{Lopez}
   F.~J.~Lopez,
   {\itshape The classification of complete minimal surfaces with
        total curvature greater than $-12\pi$},
   Trans.~Amer.~Math.~Soc. {\bfseries 334} (1992), 49--74.
\bibitem{Osserman}
    R.~Osserman,
  {\sc A Survey of Minimal Surfaces}, {2nd ed.}, Dover, 1986.
\bibitem{rs}
  W.~Rossman, K.~Sato,
  {\itshape Constant mean curvature surfaces with two ends in hyperbolic
  space}, Experimental Math., {\bfseries 7(2)} (1998), 101--119.
\bibitem{ruy1}
  W.~Rossman, M.~Umehara and K.~Yamada,
  {\itshape Irreducible constant mean curvature $1$ surfaces in
    hyperbolic space with positive genus}, 
  T\^ohoku Math.~J. {\bfseries 49}  (1997), 449--484.
\bibitem{ruy2}
  \bysame,
  {\itshape A new flux for mean curvature $1$ surfaces
    in hyperbolic $3$-space, and applications},
    Proc.~Amer.~Math.~Soc. {\bfseries 127} (1999), 2147--2154.
\bibitem{ruy3}
  \bysame,
  {\itshape Mean curvature $1$ surfaces in hyperbolic
   $3$-space with low total curvature I},
   preprint, math.DG/0008015.
\bibitem{ruy4}
  \bysame,
  {\itshape Period problems for mean curvature one surfaces
   in $H^3$ {\rm (}with application to surfaces of low total 
   curvature{\rm )}}, preprint.  
\bibitem{Sm} 
  A.~J.~Small,
  {\itshape Surfaces of Constant Mean Curvature 
       $1$ in $H^3$ and Algebraic Curves on a Quadric}, 
  Proc.~Amer.~Math.~Soc. 
  {\bfseries 122} (1994), 1211--1220.
\bibitem{Troyanov1}
  M.~Troyanov,
  {\itshape Metric of constant curvature on a sphere with two conical
       singularities}, 
  in ``Differential Geometry'', Lect. Notes in Math. vol.~1410, 
  Springer-Verlag,  (1989), 296--306. 
\bibitem{Troyanov2}
   \bysame,
  {\itshape Prescribing curvature on compact surfaces with conical
          singularities},
   Trans. Amer.~Math.~Soc. {\bfseries 324} (1991), 793--821.
\bibitem{uy1}
  M.~Umehara and K.~Yamada,
  {\itshape Complete surfaces of constant mean curvature-$1$
       in the hyperbolic $3$-space},
  {Ann. of Math. {\bfseries 137} (1993), 611--638.}
\bibitem{uy2}
  \bysame,
  {\itshape A parameterization of Weierstrass formulae and 
     perturbation of some complete minimal surfaces of
        $\R^3$ into the hyperbolic $3$-space},
   J. reine u.~angew.~Math. {\bfseries 432} (1992), 93--116.
\bibitem{uy3}
   \bysame,
   {\itshape Surfaces of constant mean curvature-$c$
        in $H^3(-c^2)$ with prescribed hyperbolic Gauss map},
   Math. Ann. {\bfseries 304} (1996), 203--224.
\bibitem{uy4}
   \bysame,
   {\itshape Another construction of a CMC-$1$ surface in $H^3$},
    Kyungpook Math. J. {\bfseries 35} (1996), 831--849.
\bibitem{uy5}
   \bysame,
   {\itshape A duality on CMC-$1$ surface in the hyperbolic $3$-space
        and a hyperbolic analogue of the Osserman Inequality},
    Tsukuba J. Math. {\bfseries 21} (1997), 229--237.
\bibitem{uy6}
   \bysame,
   {\itshape Metrics of constant curvature one with
        three conical singularities on the $2$-sphere}, 
   Illinois J. Math.~{\bfseries 44} (2000), 72--94.
\bibitem{Yu}
    Z.~Yu,
   {\itshape Value distribution of hyperbolic Gauss maps},
  Proc.~Amer.~Math.~Soc. 
  {\bfseries 125} (1997), 2997--3001.
\bibitem{Yu2}
   \bysame,
   {\itshape The inverse surface and the Osserman Inequality},
  Tsukuba J. Math. 
  {\bfseries 22} (1998), 575--588.
\end{thebibliography}
\end{document}